# HIGH-POWER ASYMPTOTICS OF SOME WEIGHTED HARMONIC BERGMAN KERNELS

Miroslav Engliš

ABSTRACT. For weights $\rho$ which are either radial on the unit ball or depend only on the vertical coordinate on the upper half-space, we describe the asymptotic behaviour of the corresponding weighted harmonic Bergman kernels with respect to $\rho^\alpha$ as $\alpha \to +\infty$. This can be compared to the analogous situation for the holomorphic case, which is of importance in the Berezin quantization as well as in complex geometry.

## 1. INTRODUCTION

Let $\Omega$ be a domain in $\mathbf{C}^n$, $\rho$ a positive smooth ($= C^\infty$) weight on $\Omega$, $L^2_{\text{hol}}(\Omega, \rho^\alpha)$ the subspace of all holomorphic functions in the weighted Lebesgue space $L^2(\Omega, \rho^\alpha)$, and $K_\alpha(x, y)$ the reproducing kernel for $L^2_{\text{hol}}(\Omega, \rho^\alpha)$, i.e. the weighted Bergman kernel on $\Omega$ with respect to the weight $\rho^\alpha$. Under suitable hypothesis on $\Omega$ and $\rho$ (namely, for $\Omega$ bounded and pseudoconvex, $\log \frac{1}{\rho}$ strictly plurisubharmonic, and $\rho$ a defining function for $\Omega$, i.e. vanishing to precisely the first order at the boundary), it is then known that

$$(1) \qquad K_\alpha(x,x) \sim \frac{\alpha^n}{\pi^n \rho(x)^\alpha} \det\left[\partial\overline{\partial}\log\frac{1}{\rho(x)}\right] \qquad \text{as } \alpha \nearrow +\infty.$$

In fact, there is even a similar result for $K_\alpha(x, y)$ with $y$ close to $x$, and one also has a complete asymptotic expansion as $\alpha \nearrow +\infty$

$$(2) \qquad K_\alpha(x,y) \approx \frac{\alpha^n}{\pi^n \rho(x,y)^\alpha} \sum_{j=0}^{\infty} \frac{b_j(x,y)}{\alpha^j}, \qquad b_0(x,x) = \det[\partial\overline{\partial}\log\tfrac{1}{\rho(x)}],$$

with some "sesqui-analytic extension" $\rho(x, y)$ of $\rho(x)$ and sesqui-analytic coefficient functions $b_j(x, y)$. Furthermore one can differentiate (1) and (2) termwise any number of times. There are, finally, variants also for the weighted Bergman spaces with respect to $\rho^\alpha \psi^m$, where $\psi$ is another weight function satisfying the same hypotheses as $\rho$ and $m \geq 0$ is a fixed real number. All these "high power asymptotics" can also be extended from functions on domains $\Omega$ to sections of holomorphic Hermitian line bundles over a manifold $\Omega$, and are then of central importance in certain approaches to quantization (the Berezin-Toeplitz quantization procedure), as well as in complex geometry (where (1) is sometimes known as the Tian-Yau-Zelditch

1991 *Mathematics Subject Classification.* Primary 46E22; Secondary 31B05, 32A36.
Research supported by GA ČR grant no. 201/12/0426 and RVO funding for IČ 67985840.

Typeset by $\mathcal{A}_{\mathcal{M}}\mathcal{S}$-TEX

expansion, and plays prominent role e.g. in connection with semistability and constant scalar curvature metrics on $\Omega$); see for instance Berezin [2], Engliš [7], [9], Zelditch [16], Catlin [4], Donaldson [6], and the references therein.

While there exist several well-understood variants of methods how to prove (1) (or (2)) nowadays, none of them makes it quite clear what does the holomorphy of functions in $L^2_{\text{hol}}$ have to do with (1), (2) or with the coefficients $b_j$ above; in fact, a priori there is little reason to expect that holomorphic functions should have anything to do either with quantization or with constant scalar curvature metrics, and one is just left to wonder at Berezin's original insight in noticing (1) and its applications. In particular, it remains quite elusive what happens for other reproducing kernel subspaces in $L^2(\Omega, \rho^\alpha)$.

The goal of this paper is to explore the analogue of (1) for the spaces of harmonic, rather than holomorphic, functions, i.e. for the reproducing kernels $R_\alpha(x, y)$ — the harmonic Bergman kernels — of the subspaces $L^2_{\text{harm}}(\Omega, \rho^\alpha)$ of all harmonic functions in $L^2(\Omega, \rho^\alpha)$.

In the holomorphic setting, the simplest examples for (1) and (2) are the standard weighted Bergman spaces on the unit disc $\mathbf{D}$ in $\mathbf{C}$ with $\rho(z) = 1 - |z|^2$, when

$$(3) \qquad K_\alpha(x, y) = \tfrac{\alpha+1}{\pi}(1 - x\overline{y})^{-\alpha-2};$$

or, equivalently (via the Cayley transform), on the upper half-plane $\{z : \operatorname{Im} z > 0\}$ in $\mathbf{C}$ with $\rho(z) = \operatorname{Im} z$ and

$$(4) \qquad K_\alpha(x, y) = \frac{\alpha+1}{4\pi} \left(\frac{x - \overline{y}}{2i}\right)^{-\alpha-2}.$$

More generally, for the unit ball $\mathbf{B}^{2n}$ of $\mathbf{C}^n \cong \mathbf{R}^{2n}$ with $\rho(z) = 1 - |z|^2$ one gets

$$(5) \qquad K_\alpha(x, y) = \frac{\Gamma(\alpha + n + 1)}{\Gamma(\alpha + 1)\pi^n}(1 - x\overline{y})^{-\alpha-n-1}.$$

Explicit formulas for the harmonic analogues of (3)–(5), namely for $R_\alpha(x, y)$ for the upper half-space $\mathbf{H}^n = \{x \in \mathbf{R}^n : x_n > 0\}$ with the weight $\rho(x) = x_n$, and for the unit ball $\mathbf{B}^n = \{x \in \mathbf{R}^n : |x| < 1\}$ with $\rho(x) = 1 - |x|^2$, have been computed in many places, see e.g. Coifman and Rochberg [5], Jevtic and Pavlovic [14], Miao [15], or the book by Axler, Bourdon and Ramey [1]. For $\Omega = \mathbf{B}^n$ and $\alpha = 0$ (i.e. the unweighted situation), the kernel is given by

$$R_0(x, y) = \frac{\Gamma(\tfrac{n}{2})}{2\pi^{n/2}} \frac{(n-4)|x|^4|y|^4 + (8\langle x, y\rangle - 2n - 4)|x|^2|y|^2 + n}{(1 - 2\langle x, y\rangle + |x|^2|y|^2)^{n/2+1}}.$$

For the weighted case with $\alpha > -1$, one already gets the much more complicated formula

$$R_\alpha(x, y) = \frac{\Gamma(\alpha + \tfrac{n}{2} + 1)}{\Gamma(\alpha + 1)\pi^{n/2}} F_1\left(\begin{matrix} \alpha + \tfrac{n}{2} + 1 \\ \tfrac{n}{2} - 1\end{matrix}; \tfrac{n}{2} - 1, \tfrac{n}{2} - 1; z, \overline{z}\right),$$

involving Appel's hypergeometric function $F_1$ [3]; here $z = x \cdot y + i\sqrt{|x|^2|y|^2 - (x \cdot y)^2}$. For $x = y$, this reduces to the ordinary hypergeometric function

$$R_\alpha(x, x) = \frac{\Gamma(\alpha + \tfrac{n}{2} + 1)}{\Gamma(\alpha + 1)\pi^{n/2}} {}_2F_1\left(\begin{matrix}\alpha + \tfrac{n}{2} + 1, n - 2 \\ \tfrac{n}{2} - 1\end{matrix}\bigg| |x|^2\right)$$

from which one gets the asymptotics

(6) $$R_\alpha(x,x) \sim \begin{cases} \dfrac{2\Gamma(\frac{n}{2})}{\pi^{n/2}\Gamma(n-1)} \dfrac{\alpha^{n-1}|x|^{n-2}}{(1-|x|^2)^{n+\alpha}} & \text{for } x \neq 0 \\ \dfrac{\alpha^{n/2}}{\pi^{n/2}} & \text{for } x = 0 \end{cases} \quad \text{as } \alpha \nearrow +\infty,$$

as the simplest harmonic analogue of (1). Similarly, for the upper half-space and $\alpha = 0$, the unweighted kernel is given by

$$R_0(x,y) = \frac{2\Gamma(\frac{n}{2})}{\pi^{n/2}} \frac{(n-1)(x_n+y_n)^2 + (x_n-y_n)^2 - |x-y|^2}{[(x_n+y_n)^2 - (x_n-y_n)^2 + |x-y|^2]^{n/2+1}},$$

while for general $\alpha$ one can compute e.g. from [13] that

(7) $$R_\alpha(x,x) = \frac{\Gamma(n+\alpha)2^{3-2n}}{\pi^{\frac{n-1}{2}}\Gamma(\alpha+1)\Gamma(\frac{n-1}{2})x_n^{n+\alpha}} \sim \frac{\alpha^{n-1}2^{3-2n}}{\pi^{\frac{n-1}{2}}\Gamma(\frac{n-1}{2})x_n^{n+\alpha}},$$

giving the harmonic analogue of (1) for the upper half-space with $\rho(x) = x_n$. Finally, one can also consider the entire $\mathbf{R}^n$ with the Gaussian weight $\rho(x) = e^{-|x|^2}$ (the harmonic Fock, or Segal-Bargmann, space), in which case it was derived in [11] that

$$R_\alpha(x,y) = \frac{\alpha^{n/2}}{\pi^{n/2}} \Phi_2\left(\begin{array}{c}\frac{n}{2}-1, \frac{n}{2}-1 \\ \frac{n}{2}-1\end{array}\bigg| \alpha z, \alpha \overline{z}\right),$$

with Horn's hypergeometric function $\Phi_2$ and again $z = x \cdot y + i\sqrt{|x|^2|y|^2 - (x \cdot y)^2}$; for $x = y$ this reduces to the confluent hypergeometric function

$$R_\alpha(x,x) = \frac{\alpha^{n/2}}{\pi^{n/2}} {}_1F_1\left(\begin{array}{c}n-2 \\ \frac{n}{2}-1\end{array}\bigg| \alpha|x|^2\right),$$

yielding

(8) $$R_\alpha(x,x) \sim \begin{cases} \dfrac{2\Gamma(\frac{n}{2})}{\Gamma(n-1)\pi^{n/2}} e^{\alpha|x|^2}|x|^{n-2}\alpha^{n-1} & \text{for } x \neq 0, \\ \dfrac{\alpha^{n/2}}{\pi^{n/2}} & \text{for } x = 0, \end{cases}$$

as $\alpha \nearrow +\infty$. Note that in (6) and (8), we get the "Stokes phenomenon" of different asymptotics at $x = 0$ and $x \neq 0$, which is unparalleled in the holomorphic case as well as in the case of the upper half-space in (7).

Our result here is a rather coarse description for the asymptotics of $R_\alpha(x,x)$ for fairly general $\rho$ and $\Omega$, and more precise descriptions on the level of (1) for domains and weights of a particular form.

**Theorem 1.** *Let $\Omega \subset \mathbf{R}^n$, $n \geq 2$, be a bounded domain and $\rho$ a bounded positive continuous function on $\Omega$ such that $\log \frac{1}{\rho}$ is convex. Then*

$$\lim_{\alpha \to \infty} R_\alpha(x,x)^{1/\alpha} = \frac{1}{\rho(x)}.$$

Keeping the usual definition from complex analysis, we call $\rho$ a *defining function* for $\Omega$ if $\rho > 0$ on $\Omega$ and $\rho$ vanishes precisely to the first order at the boundary $\partial\Omega$, i.e. $\rho = 0 < \|\nabla\rho\|$ on $\partial\Omega$.

**Theorem 2.** *Let $\Omega = \mathbf{B}^n$, $n \geq 2$, and let $\rho$ be radial, i.e. $\rho(x) = \phi(|x|^2)$ for some positive $\phi \in C^\infty[0,1]$. Assume that $\rho$ is a defining function (i.e. $\phi(1) = 0$ an $\phi'(1) < 0$) and that $(\frac{t\phi'}{\phi})' < 0$. Then for any $x \neq 0$,*

$$(9) \quad \lim_{\alpha \to \infty} \alpha^{1-n} \rho(x)^\alpha R_\alpha(x,x) = \frac{2\Gamma(\frac{n}{2}) t^{\frac{n}{2}-1}}{\pi^{n/2} \Gamma(n-1)} \Big(-\frac{\phi'}{\phi}\Big)^{n-2} \Big(-\frac{t\phi'}{\phi}\Big)' \Big|_{t=|x|^2}.$$

Let us call a positive function $g$ on $(0, +\infty)$ admissible if $\int_0^\infty e^{tx} g(x)\, dx = +\infty$ for all $t > 0$; this means that $g$ should not decay too rapidly at infinity.

**Theorem 3.** *Let $\Omega = \mathbf{B}^n$, $n \geq 2$, and assume that $\rho(x) = \rho(x_n)$ depends only on the vertical coordinate, is admissible, vanishes at $x_n = 0$ precisely to the first order, and $\rho' > 0$, $(\rho'/\rho)' < 0$ on $(0, +\infty)$. Then*

$$(10) \quad \lim_{\alpha \to \infty} \alpha^{1-n} \rho(x)^\alpha R_\alpha(x,x) = \frac{2^{3-2n}}{\pi^{\frac{n-1}{2}} \Gamma(\frac{n-1}{2})} \Big(\frac{\rho'}{\rho}\Big)^{n-2} \Big(-\frac{\rho'}{\rho}\Big)'.$$

Note that the choices $\phi(t) = 1-t$ and $\rho(x) = x_n$ recover (6) and (8), respectively. In fact, $\phi(t) = e^{-t}$ recovers also (7).

The proof of Theorem 1 appears in Section 2, and those of Theorems 2 and 3 in Sections 3 and 4, respectively. Our main idea is to reduce the harmonic case to the holomorphic one (on a different domain) and then use (1). Some concluding remarks are given in the final Section 5.

Throughout the rest of the paper, we write just $K_\alpha(x)$, $R_\alpha(x)$ for $K_\alpha(x,x)$ and $R_\alpha(x,x)$, respectively; and, as usual, "$A(x) \sim B(x)$ as $\alpha \to +\infty$" means that $A(x)/B(x) \to 1$ as $\alpha \to +\infty$. The norm in $L^2(\Omega, \rho^\alpha)$ is denoted by $\|\cdot\|_\alpha$, and "plurisubharmonic" will be abbreviated to "psh".

## 2. Coarse asymptotics

The proof of our first theorem is actually almost the same as for the holomorphic case in [8].

*Proof of Theorem 1.* Let $D(x,r)$ be the polydisc with center $x$ and radius $r$, where $r > 0$ is so small that $D(x,r) \subset \Omega$. By the mean value property, for any $h \in L^2_{\text{harm}}(\Omega, \rho^\alpha)$,

$$h(x) = (\pi r^2)^{-n} \int_{D(x,r)} h(y)\, dy,$$

where $dy$ stands for the Lebesgue measure. The Cauchy-Schwarz inequality gives

$$|h(x)| \leq (\pi r^2)^{-n} \Big(\int_{D(x,r)} |h|^2 \rho^\alpha\, dy\Big)^{1/2} \Big(\int_{D(x,r)} \rho^{-\alpha}\, dy\Big)^{1/2}$$

$$\leq (\pi r^2)^{-n/2} \|h\|_\alpha \big(\sup_{D(x,r)} \tfrac{1}{\rho}\big)^{\alpha/2}.$$

Now by the extremal property of reproducing kernels, $R_\alpha(x)^{1/2}$ is the norm of the evaluation functional $f \mapsto f(x)$ on $L^2_{\text{harm}}(\Omega, \rho^\alpha)$, that is,

$$(11) \quad R_\alpha(x) = \sup\{|h(x)|^2 : h \in L^2_{\text{harm}}(\Omega, \rho^\alpha), \|h\| \leq 1\}.$$

Thus
$$R_\alpha(x) \leq (\pi r^2)^{-n} (\sup_{D(x,r)} \tfrac{1}{\rho})^\alpha.$$

Taking $\alpha$-th roots on both sides and letting $\alpha \nearrow +\infty$ gives
$$\limsup_{\alpha \to \infty} R_\alpha(x)^{1/\alpha} \leq \sup_{D(x,r)} \tfrac{1}{\rho}.$$

Letting $r \searrow 0$, the continuity of $\rho$ implies

(12) $$\limsup_{\alpha \to \infty} R_\alpha(x)^{1/\alpha} \leq \frac{1}{\rho(x)}.$$

On the other hand, by (11), for any $h \in L^2_{\text{harm}}(\Omega, \rho^\alpha)$ not identically zero we have
$$R_\alpha(x) \geq \frac{|h(x)|^2}{\|h\|_\alpha^2}.$$

Take in particular $h(x) = e^{\alpha(x \cdot z + c)}$, with arbitrary $c \in \mathbf{R}$ and $z \in \mathbf{C}^n$ satisfying $z_1^2 + z_2^2 + \cdots + z_n^2 =: z \cdot z = 0$; clearly this is a bounded function and hence in $L^2(\Omega, \rho^\alpha)$ (since $\Omega$ and $\rho$ are bounded by hypothesis), while the condition $z \cdot z = 0$ ensures that $h$ is harmonic. Thus
$$R_\alpha(x) \geq \frac{e^{2\alpha(x \cdot \operatorname{Re} z + c)}}{\| |e^{x \cdot z + c}|^2 \rho \|_{L^\alpha(\Omega)}^\alpha}$$

and
$$R_\alpha(x)^{1/\alpha} \geq \frac{e^{2(x \cdot \operatorname{Re} z + c)}}{\| |e^{x \cdot z + c}|^2 \rho \|_{L^\alpha(\Omega)}}.$$

Now since $|e^{x \cdot z + c}|^2 \rho(x)$ is bounded and $\Omega$ has finite Lebesgue measure, it is standard that $\| |e^{x \cdot z + c}|^2 \rho \|_{L^\alpha(\Omega)} \to \| |e^{x \cdot z + c}|^2 \rho \|_{L^\infty(\Omega)}$ as $\alpha \nearrow +\infty$. We thus obtain

(13) $$\liminf_{\alpha \to +\infty} R_\alpha(x)^{1/\alpha} \geq \sup\{e^{2(x \cdot \operatorname{Re} z + c)} : z \cdot z = 0, c \in \mathbf{R}, |e^{\bullet \cdot z + c}|^2 \rho \leq 1\}.$$

Writing $z = a + bi$ for the real and imaginary parts of $z$, the condition $z \cdot z = 0$ becomes

(14) $$a \cdot a = b \cdot b, \qquad a \cdot b = 0.$$

Since $n \geq 2$, we can for any $a \in \mathbf{R}^n$ find $b \in \mathbf{R}^n$ such that (14) holds (just take any $b$ orthogonal to $a$ and of the same length). Thus (13) translates into
$$\log \liminf_{\alpha \to +\infty} R_\alpha(x)^{1/\alpha} \geq \sup\{\psi : \psi \text{ affine on } \mathbf{R}^n, \psi \leq \log \tfrac{1}{\rho}\}.$$

If $\log \tfrac{1}{\rho}$ is convex, then the right-hand side equals $1/\rho$, whence

(15) $$\liminf_{\alpha \to +\infty} R_\alpha(x)^{1/\alpha} \geq \frac{1}{\rho(x)}.$$

Combining (15) and (12), the assertion follows. □

3. The case of the ball

Let us recall some prerequisites, available e.g. in [1]. For $k = 0, 1, 2, \ldots$, let $\mathcal{H}^k(\mathbf{R}^n)$ denote the space of all harmonic polynomials on $\mathbf{R}^n$ homogeneous of degree $k$. Each such polynomial is uniquely determined by its restriction to the unit sphere $\mathbf{S}^{n-1} := \partial \mathbf{B}^n \subset \mathbf{R}^n$, and we denote by $\mathcal{H}^k$ the space of all such restrictions (called "spherical harmonics of degree $k$"), viewed as a subspace of $L^2(\mathbf{S}^{n-1}, d\sigma)$, where $d\sigma$ stands for the normalized surface measure on $\mathbf{S}^{n-1}$. Then $\mathcal{H}^k \perp \mathcal{H}^l$ for $k \neq l$, and the span of all $\mathcal{H}^k$, $k \geq 0$, is dense in $L^2(\mathbf{S}^{n-1}, d\sigma)$. Furthermore, each $\mathcal{H}^k$ is a reproducing kernel space, with reproducing kernel $Z_k(x, y)$ given by the so-called zonal harmonic; this is a certain Gegenbauer polynomial in $x \cdot y$. Finally, any harmonic function $f$ on $\mathbf{B}^n$ can be uniquely expressed in the form

$$f = \sum_{k=0}^{\infty} f_k, \qquad f_k \in \mathcal{H}^k(\mathbf{R}^n),$$

with the sum converging uniformly on compact subsets of $\mathbf{B}^n$.

For any radial weight $\rho(x) = \phi(|x|^2)$ on $\mathbf{B}^n$ and $0 < r < 1$, we therefore get (recalling that the volume of $\mathbf{S}^{n-1}$ equals $2\pi^{n/2}/\Gamma(\frac{n}{2})$)

$$\begin{aligned}
\int_{r\mathbf{B}^n} |f|^2 \rho \, dx &= \sum_{j,k} \int_{r\mathbf{B}^n} f_j \overline{f_k} \rho \, dx \\
&= \sum_{j,k} \frac{2\pi^{n/2}}{\Gamma(\frac{n}{2})} \int_0^r r^{j+k} \phi(r^2) \, r^{n-1} \int_{\mathbf{S}^{n-1}} f_j(\zeta) \overline{f_k(\zeta)} \, d\sigma(\zeta) \, dr \\
&= \frac{2\pi^{n/2}}{\Gamma(\frac{n}{2})} \sum_k \Big( \int_0^r r^{2k+n-1} \phi(r^2) \, dr \Big) \|f_k\|^2_{L^2(\mathbf{S}^{n-1}, d\sigma)},
\end{aligned}$$

by the orthogonality of $\mathcal{H}^k$ and $\mathcal{H}^l$ for $k \neq l$. Denoting

$$(16) \qquad \int_0^1 r^k \phi(r^2) \, dr =: \rho_k,$$

it follows upon letting $r \nearrow 1$ that

$$\|f\|^2_{L^2_{\text{harm}}(\Omega, \rho)} = \frac{2\pi^{n/2}}{\Gamma(\frac{n}{2})} \sum_k \rho_{2k+n-1} \|f_k\|^2_{L^2_{\text{harm}}(\Omega, \rho)},$$

and, consequently, the reproducing kernel of $L^2_{\text{harm}}(\Omega, \rho)$ is given by

$$R(x, y) = \frac{\Gamma(\frac{n}{2})}{2\pi^{n/2}} \sum_k \frac{Z_k(x, y)}{\rho_{2k+n-1}}.$$

For $x = y$, it is known that $Z_k(x, x) = N_{k,n} |x|^{2k}$, where $N_{k,n} = \dim \mathcal{H}^k$ is given by

$$N_{k,n} = \frac{(n+k-3)!(n+2k-2)}{k!(n-2)!}.$$

We thus obtain

$$R(x) = \frac{\Gamma(\frac{n}{2})}{2\pi^{n/2}} \sum_k \frac{N_{k,n}}{\rho_{2k+n-1}} |x|^{2k}. \tag{17}$$

On the other hand, a completely similar formula is available also for the holomorphic Bergman kernel $K(x,y)$ on the unit ball $\mathbf{B}^{2m}$ of $\mathbf{R}^{2m} \cong \mathbf{C}^m$, $m \geq 1$, with respect to a radial weight function $w(z) = \psi(|z|^2)$. Namely, it is standard that the monomials $z^\nu$, $\nu$ a multiindex, are then orthogonal, with norm squares

$$\int_{\mathbf{B}^{2n} \subset \mathbf{C}^n} |z^\nu|^2 \, w(z) \, dz = \underbrace{\left(\int_0^1 r^{2|\nu|+2m-1} \psi(r^2) \, dr\right)}_{=:w_{2|\nu|+2m-1}} \left(\int_{\mathbf{S}^{2n-1}} |\zeta^\nu|^2 \, d\zeta\right)$$

$$= \frac{2\pi^m \nu!}{(|\nu|+m-1)!} w_{2|\nu|+2m-1},$$

and by the familiar formula expressing the reproducing kernel in terms of an arbitrary orthonormal basis,

$$K(x,y) = \sum_\nu \frac{x^\nu \overline{y}^\nu}{2\pi^m \nu!} \frac{(|\nu|+m-1)!}{w_{2|\nu|+2m-1}} = \frac{1}{2\pi^m} \sum_{k=0}^\infty \frac{\langle x,y \rangle^k}{k!} \frac{(k+m-1)!}{w_{2k+2m-1}}.$$

We thus obtain

$$K(z) = \frac{1}{2\pi^m} \sum_{k=0}^\infty \frac{(k+m-1)!}{k! w_{2k+2m-1}} |z|^{2k}. \tag{18}$$

Finally, note that for any function

$$F(t) = \sum_{k=0}^\infty c_k t^k$$

holomorphic on the unit disc, we have

$$\sum_{k=0}^\infty \frac{(k+m-1)!}{k!} c_k t^k = (t^{m-1} F)^{(m-1)}, \tag{19}$$

and

$$\sum_{k=0}^\infty N_{k,n} c_k t^k = \frac{(t^{n-2} F)^{(n-2)} + t(t^{n-3} F)^{(n-2)}}{(n-2)!}, \tag{20}$$

as can be checked by elementary manipulations.

The last thing we will need is the fact that the condition $(\frac{t\phi'}{\phi})' < 0$ is actually equivalent to the function $\log \frac{1}{\phi(|z|^2)}$ being strictly-psh on $\mathbf{B}^{2m}$, see e.g. [7], Section 3. Since $\phi \in C^\infty[0,1]$ by hypothesis, this condition also implies that $-\frac{t\phi'}{\phi} > 0$ for $t > 0$, and hence $\phi' < 0$ for $t > 0$, that is, $\phi$ is decreasing on $(0,1)$.

*Proof of Theorem 2.* Assume first that $n$ is even. Take $m = \frac{n}{2}$ in (18), with $\psi = \phi$. Replacing $\phi$ by $\phi^\alpha$ and denoting the corresponding $\rho_k$ from (16) by $\rho_k(\alpha)$, we thus get

$$R_\alpha(x) = \frac{\Gamma(m)}{2\pi^m} \sum_{k=0}^{\infty} \frac{N_{k,2m}}{\rho_{2k+2m-1}(\alpha)} t^k, \qquad K_\alpha(z) = \frac{1}{2\pi^m} \sum_{k=0}^{\infty} \frac{(k+m-1)!}{k!\rho_{2k+2m-1}(\alpha)} t^k,$$

where we have set $t = |x|^2$ and $t = |z|^2$, respectively. Thus $\frac{2\pi^m}{\Gamma(m)} R_\alpha(x) =: r_\alpha(t)$ and $(2\pi^m) K_\alpha(z) =: k_\alpha(t)$ are related as in (19) and (20), with

$$F(t) \equiv F_\alpha(t) = \sum_{k=0}^{\infty} \frac{t^k}{\rho_{2k+2m-1}(\alpha)}.$$

By the Leibniz rule, $r_\alpha(t)$ comes as a sum of terms of the form $c_j t^{m_j} k_\alpha^{(j)}(t)$, $j = -(m-1), \ldots, m-1$, with some real numbers $c_j$ and integers $m_j$ (independent of $\alpha$), where $k_\alpha^{(j)}$ denotes the $j$-th derivatives of $k_\alpha$ for $j \geq 0$ and the $|j|$-th primitive of $k_\alpha$ (normalized to vanish to order $|j|$ at $t = 0$) for $j < 0$. Now by (1),

$$(21) \qquad k_\alpha(|z|^2) \sim \frac{2\alpha^m}{\rho(z)^\alpha} \det\left[\partial\overline{\partial} \log \frac{1}{\rho(z)}\right],$$

and this also remains in force upon applying any derivative to both sides. Since a short computation shows that

$$\det\left[\partial\overline{\partial} \log \frac{1}{\rho(z)}\right] = \left(-\frac{\phi'}{\phi}\right)^{m-1} \left(-\frac{t\phi'}{\phi}\right)' \bigg|_{t=|z|^2},$$

we have

$$k_\alpha \sim \frac{2\alpha^m}{\phi^\alpha} \left(-\frac{\phi'}{\phi}\right)^{m-1} \left(-\frac{t\phi'}{\phi}\right)'.$$

Differentiation gives

$$k'_\alpha \sim \frac{2\alpha^m}{\phi^\alpha} \frac{-\alpha\phi'}{\phi} \left(-\frac{\phi'}{\phi}\right)^{m-1} \left(-\frac{t\phi'}{\phi}\right)' + \frac{2\alpha^m}{\phi^\alpha} \left[\left(-\frac{\phi'}{\phi}\right)^{m-1} \left(-\frac{t\phi'}{\phi}\right)'\right]'.$$

As $\alpha \nearrow +\infty$, the first term dominates the second. Thus by induction

$$(22) \qquad k_\alpha^{(j)} \sim \left(\frac{-\alpha\phi'}{\phi}\right)^j k_\alpha$$

for any $j \geq 0$. On the other hand,

$$k_\alpha^{(-1)}(t) = \int_0^t k_\alpha \sim 2\alpha^m \int_0^t \frac{1}{\phi^\alpha} \left(-\frac{\phi'}{\phi}\right)^{m-1} \left(-\frac{t\phi'}{\phi}\right)'.$$

The right-hand side is a standard Laplace-type integral, that is, an integral of the form

$$(23) \qquad I(\alpha) = \int_a^b F(x) e^{\alpha S(x)} \, dx$$

with real-valued function $S$. It is known that $I(\alpha)$ gets the largest contribution from points where $S$ attains its maximum, and, in particular, if the maximum is attained at the endpoint $x = b$ and $S'(b) > 0$, then (see e.g. [12], §II.1.4)

$$I(\alpha) \sim \frac{F(b)}{S'(b)} \frac{e^{\alpha S(b)}}{\alpha}.$$

Since, as we have observed, $\phi$ is decreasing, $\frac{1}{\phi}$ indeed attains its maximum at the endpoint $t$, and thus

$$k_\alpha^{(-1)} \sim \frac{2\alpha^m}{\phi^\alpha} \Big(-\frac{\phi'}{\phi}\Big)^{m-1} \Big(-\frac{t\phi'}{\phi}\Big)' \cdot \Big(\frac{\phi}{-\alpha\phi'}\Big).$$

Proceeding inductively, it follows that (22) in fact remains in force for $j \leq 0$ as well. Thus the leading term in the asymptotics of $r_\alpha$ as $\alpha \nearrow +\infty$ will be the one coming from (22) with $j = m-1$. In other words, $k_\alpha \sim t^{m-1} F^{(m-1)}$ and

$$r_\alpha \sim \frac{2t^{2m-2} F_\alpha^{(2m-2)}}{(2m-2)!} \sim \frac{2t^{m-1} k_\alpha^{(m-1)}}{(2m-2)!} \sim \frac{2t^{m-1}}{(2m-2)!} \Big(\frac{-\alpha\phi'}{\phi}\Big)^{m-1} k_\alpha$$
$$\sim \frac{4\alpha^{2m-1} t^{m-1}}{(2m-2)! \phi^\alpha} \Big(-\frac{\phi'}{\phi}\Big)^{2m-2} \Big(-\frac{t\phi'}{\phi}\Big)',$$

or

$$R_\alpha(x) = \frac{\Gamma(m)}{2\pi^m} r_\alpha(t) \sim \frac{2\alpha^{n-1} t^{\frac{n}{2}-1} \Gamma(\frac{n}{2})}{\pi^{n/2}(n-2)! \phi^\alpha} \Big(-\frac{\phi'}{\phi}\Big)^{n-2} \Big(-\frac{t\phi'}{\phi}\Big)' \Big|_{t=|x|^2},$$

proving (9) for even $n$.

For $n$ odd, take $m = \frac{n-1}{2}$ and replace $\phi, \psi$ by $\phi^\alpha$ and $\phi^\alpha \sqrt{t}$, respectively. Then $\phi^\alpha \sqrt{t} = \phi^{\alpha-1}(\phi\sqrt{t})$ and $\phi\sqrt{t}$ again vanishes to exactly the first order at $\partial \mathbf{B}^{2m}$ and $\log \frac{1}{\phi(|z|^2)|z|}$ is still strictly-psh, except that $\phi(|z|^2)|z|$ now fails to be smooth at the origin; however, it is known that (1) — or, more precisely, its variant for two weight functions mentioned at the beginning of this paper — then still remains in force at points where the weights are smooth, so in our case for $z \neq 0$. (See [10], Theorem 1.) We thus have

$$K_{\phi(|z|^2)^\alpha |z|}(z) \sim \frac{\alpha^m}{\pi^m \phi(|z|^2)^\alpha |z|} \det\Big[\partial\overline{\partial} \log \frac{1}{\phi(|z|^2)|z|}\Big]$$

for $z \neq 0$. Arguing as in the case of even $n$, now with $k_\alpha(t) = 2\pi^m K_{\phi(|z|^2)^\alpha |z|}(z)$, $r_\alpha(t) = \frac{2\pi^{n/2}}{\Gamma(\frac{n}{2})} R_\alpha(x)$, and $F_\alpha(t) = \sum_k \frac{t^k}{\rho_{k+2m}(\alpha)}$, we arrive at

$$r_\alpha \sim \frac{2t^{2m-1} F^{(2m-1)}}{(2m-1)!} \sim \frac{2t^m k_\alpha^{(m)}}{(2m-1)!} \sim \Big(\frac{-\phi'}{\phi}\Big)^m \frac{2\alpha^m t^m}{(2m-1)!} k_\alpha$$
$$\sim \frac{4\alpha^{2m} t^m}{\phi^\alpha t^{1/2}(2m-1)!} \Big(-\frac{\phi'}{\phi}\Big)^{2m-1} \Big(-\frac{t\phi'}{\phi}\Big)'$$

and

$$R_\alpha(x) \sim \frac{2\Gamma(\frac{n}{2})\alpha^{n-1} t^{\frac{n}{2}-1}}{\pi^{n/2} \phi^\alpha (n-2)!} \Big(-\frac{\phi'}{\phi}\Big)^{n-2} \Big(-\frac{t\phi'}{\phi}\Big)',$$

which settles (9) also for odd $n$ and thus completes the proof of the theorem. $\square$

## 4. The upper half-space

We again begin by reviewing some standard prerequisites on harmonic functions on $\mathbf{H}^n$. Write points in $\mathbf{H}^n$ temporarily as $(x,y)$, with $x \in \mathbf{R}^{n-1}$ and $y > 0$, and let

$$\hat{f}_y(\xi) = \int_{\mathbf{R}^{n-1}} f(x,y) e^{-ix \cdot \xi} \, d\xi$$

denote the Fourier transform of a function $f(x,y) \equiv f_y(x)$ on $\mathbf{H}^n$. The condition that $f$ be harmonic then translates into $\frac{\partial^2}{\partial y^2}\hat{f}_y + |\xi|^2 \hat{f}_y = 0$, or

$$\hat{f}_y(\xi) = A(\xi) e^{-|\xi|y} + B(\xi) e^{|\xi|y}$$

for some functions $A, B$. Now for any weight $\rho(x,y) = \rho(y)$ depending only on the vertical coordinate $y$, we have by Parseval

$$\iint_{\mathbf{H}^n} |f|^2 \rho \, dx \, dy = \int_0^\infty \rho(y) \int_{\mathbf{R}^{n-1}} |f_y|^2 \, dx \, dy = (2\pi)^{1-n} \int_0^\infty \rho(y) \int_{\mathbf{R}^{n-1}} |\hat{f}_y|^2 \, d\xi \, dy.$$

Consequently, if in addition $\rho$ is admissible, this can only be finite if $B \equiv 0$. Thus for $f \in L^2_{\mathrm{harm}}(\mathbf{H}^n, \rho)$,

$$\hat{f}_y(\xi) = \hat{f}_0(\xi) e^{-y|\xi|},$$

where $\hat{f}_0 \equiv A$ has the obvious interpretation of the Fourier transform of the boundary value $f_0$ of $f$ at $y=0$, and we can continue the last computation with

$$\iint_{\mathbf{H}^n} |f|^2 \rho \, dx \, dy = (2\pi)^{1-n} \int_0^\infty \rho(y) \int_{\mathbf{R}^{n-1}} e^{-2y|\xi|} |\hat{f}_0(\xi)|^2 \, d\xi \, dy$$

$$\equiv (2\pi)^{1-n} \int_{\mathbf{R}^{n-1}} \tilde{\rho}(|\xi|) |\hat{f}_0(\xi)|^2 \, d\xi,$$

where

(24) $$\tilde{\rho}(t) := \int_0^\infty \rho(y) e^{-2ty} \, dy.$$

Comparing this with the Fourier inversion formula

$$f(a,b) = (2\pi)^{1-n} \int_{\mathbf{R}^{n-1}} \hat{f}_b(\xi) e^{ia \cdot \xi} \, d\xi = (2\pi)^{1-n} \int_{\mathbf{R}^{n-1}} \hat{f}_0(\xi) e^{ia \cdot \xi - b|\xi|} \, d\xi,$$

we see that the reproducing kernel $R(x,y;a,b) \equiv R_{a,b}(x,y)$ of $L^2_{\mathrm{harm}}(\mathbf{H}^n, \rho)$ satisfies $\tilde{\rho}(|\xi|)(R_{a,b})_0^\wedge(\xi) = e^{-b|\xi| - ia \cdot \xi}$, or

$$R(x,y;a,b) = (2\pi)^{1-n} \int_{\mathbf{R}^{n-1}} \frac{e^{i(x-a) \cdot \xi - (b+y)|\xi|}}{\tilde{\rho}(|\xi|)} \, d\xi.$$

In particular, for $(a,b) = (x,y)$,

(25) $$R(x,y) = (2\pi)^{1-n} \int_{\mathbf{R}^{n-1}} \frac{e^{-2y|\xi|}}{\tilde{\rho}(|\xi|)} \, d\xi = \frac{2^{2-n}}{\pi^{\frac{n-1}{2}} \Gamma(\frac{n-1}{2})} \int_0^\infty \frac{e^{-2yr}}{\tilde{\rho}(r)} r^{n-2} \, dr.$$

On the other hand, consider the (Siegel) domain in $\mathbf{C}^m$ given by

(26) $$\mathcal{S} := \{(z, x + yi) \in \mathbf{C}^{m-1} \times \mathbf{C} : y > |z|^2\},$$

and the (holomorphic) Bergman space on $\mathcal{S}$ with respect to a weight $\rho(z, x + yi) \equiv \rho(y - |z|^2)$ depending only on $y - |z|^2$. Writing functions on $\mathcal{S}$ as $f(z, x + yi) \equiv f_{z,y}(x)$ and letting $\hat{f}_{z,y}$ stand for the Fourier transform of $f_{z,y}$ with respect to $x$, the holomorphy of $f_{z,y}(x)$ in $x + yi$ translates into $\partial_y \hat{f}_{z,y} + \xi \hat{f}_{z,y} = 0$, or

$$\hat{f}_{z,y} = e^{-\xi y} \hat{f}_{z,0}(\xi)$$

for some function $\hat{f}_{z,0}(\xi)$ depending holomorphically on $z$. (Again, $e^{-\xi|z|^2} \hat{f}_{z,0}(\xi) = \hat{f}_{z,|z|^2}(\xi)$ can be interpreted as the Fourier transform of the boundary value $f_{z,|z|^2}$ of $f_{z,y}$ at $y = |z|^2$.) As before, we have by Plancherel

$$\iiint_{\mathcal{S}} |f|^2 \rho \, dz \, dx \, dy = \int_{\mathbf{C}^{m-1}} \int_{|z|^2}^{\infty} \int_{\mathbf{R}} |f_{z,y}(x)|^2 \rho(y - |z|^2) \, dx \, dy \, dz$$

$$= \frac{1}{2\pi} \int_{\mathbf{C}^{m-1}} \int_{|z|^2}^{\infty} \int_{\mathbf{R}} |\hat{f}_{z,y}(\xi)|^2 \rho(y - |z|^2) \, d\xi \, dy \, dz$$

$$= \frac{1}{2\pi} \int_{\mathbf{C}^{m-1}} \int_{|z|^2}^{\infty} \int_{\mathbf{R}} |\hat{f}_{z,0}(\xi)|^2 e^{-2\xi y} \rho(y - |z|^2) \, d\xi \, dy \, dz$$

$$= \frac{1}{2\pi} \int_{\mathbf{C}^{m-1}} \int_0^{\infty} \int_{\mathbf{R}} |\hat{f}_{z,0}(\xi)|^2 e^{-2\xi|z|^2 - 2\xi r} \rho(r) \, d\xi \, dr \, dz$$

$$= \frac{1}{2\pi} \int_{\mathbf{C}^{m-1}} \int_{\mathbf{R}} |\hat{f}_{z,0}(\xi)|^2 e^{-2\xi|z|^2} \tilde{\rho}(\xi) \, d\xi \, dz$$

with $\tilde{\rho}$ as in (24). If $\rho$ is in addition admissible, the last integral can be finite only if $\hat{f}_{z,0}$ is supported on $\xi > 0$, and

(27) $$\|f\|_{L^2(\mathcal{S}, \rho)}^2 = \frac{1}{2\pi} \int_0^{\infty} \int_{\mathbf{C}^{m-1}} |\hat{f}_{z,0}(\xi)|^2 e^{-2\xi|z|^2} \tilde{\rho}(\xi) \, d\xi \, dz.$$

For $f \in L^2_{\text{hol}}(\mathcal{S}, \rho)$, we thus see that the function $z \mapsto \hat{f}_{z,0}(\xi)$ belongs to the Fock space $\mathcal{F}_{2\xi} := L^2_{\text{hol}}(\mathbf{C}^{m-1}, e^{-2\xi|z|^2})$, for any $\xi > 0$, and

$$\|f\|_{L^2(\mathcal{S}, \rho)}^2 = \frac{1}{2\pi} \int_0^{\infty} \|\hat{f}_{\bullet, 0}(\xi)\|_{\mathcal{F}_{2\xi}}^2 \tilde{\rho}(\xi) \, d\xi.$$

Since the reproducing kernel of $\mathcal{F}_{2\xi}$ is known to be $(\frac{2\xi}{\pi})^{m-1} e^{2\xi \langle z, w \rangle}$, we have

$$\hat{f}_{w,0}(\xi) = \Big(\frac{2\xi}{\pi}\Big)^{m-1} \int_{\mathbf{C}^{m-1}} e^{2\xi \langle w, z \rangle} \hat{f}_{z,0}(\xi) e^{-2\xi|z|^2} \, dz,$$

whence for $a \in \mathbf{R}$ and $b > |w|^2$

$$f(w, a + ib) = f_{w,b}(a) = \frac{1}{2\pi} \int_{\mathbf{R}} \hat{f}_{w,b}(\xi) e^{ia\xi} \, d\xi = \frac{1}{2\pi} \int_{\mathbf{R}} \hat{f}_{w,0}(\xi) e^{ia\xi - b\xi} \, d\xi$$

$$= \frac{1}{2\pi} \int_0^{\infty} \Big(\frac{2\xi}{\pi}\Big)^{m-1} \int_{\mathbf{C}^{m-1}} \hat{f}_{z,0}(\xi) e^{2\xi \langle w, z \rangle - 2\xi|z|^2} e^{ia\xi - b\xi} \, dz \, d\xi.$$

Comparing this with (27), we see that the reproducing kernel $K(z, x+yi; w, a+ib) \equiv K_{w,a,b}(z,x,y)$ of $L^2_{\text{hol}}(\mathcal{S}, \rho)$ satisfies $\tilde{\rho}(\xi)(K_{w,a,b})^\wedge_{z,0}(\xi) = (\frac{2\xi}{\pi})^{m-1} e^{2\xi\langle z,w\rangle - ia\xi - b\xi}$, or

$$K(z, x+iy; w, a+ib) = \frac{1}{2\pi} \int_0^\infty \left(\frac{2\xi}{\pi}\right)^{m-1} \frac{e^{i(x-a)\xi - (b+y)\xi + 2\xi\langle z,w\rangle}}{\tilde{\rho}(\xi)} d\xi.$$

In particular, for $(w, a+ib) = (z, x+iy)$,

$$(28) \qquad K(z, x+yi) = \frac{2^{m-2}}{\pi^m} \int_0^\infty \xi^{m-1} \frac{e^{-2(y-|z|^2)\xi}}{\tilde{\rho}(\xi)} d\xi.$$

Note finally that the hypotheses

$$\rho' > 0, \quad (\rho'/\rho)' < 0$$

mean precisely that the function $\log \frac{1}{\rho(\operatorname{Im} w - |z|^2)}$ of $(z,w) \in \mathbf{C}^{m-1} \times \mathbf{C}$ is strictly-psh on $\mathcal{S}$. Indeed, denoting momentarily $\phi = \log \frac{1}{\rho}$ for brevity, the complex Hessian matrix of $\phi(\frac{w-\overline{w}}{2i} - |z|^2)$ is given by

$$\begin{bmatrix} \frac{1}{4}\phi''(t) & -\frac{z}{2i}\phi''(t) \\ \frac{\overline{z}}{2i}\phi''(t) & -\phi'(t)I + \phi''(t) z \otimes \overline{z} \end{bmatrix}, \qquad t = \frac{w-\overline{w}}{2i} - |z|^2.$$

Multiplying the first column by $\frac{2z_j}{i}$ and adding it to the $(j+1)$-st column, and similarly for rows, shows that this matrix is positive definite if and only if $\begin{bmatrix} \frac{1}{4}\phi'' & 0 \\ 0 & -\phi'I \end{bmatrix}$ is positive definite. However the latter is clearly equivalent to $\phi'' > 0$ and $\phi' < 0$, establishing the claim.

The last argument also shows that

$$(29) \qquad \det\left[\partial\overline{\partial}\phi(\operatorname{Im} w - |z|^2)\right] = \frac{\phi''(t)}{4}(-\phi'(t))^{m-1}\Big|_{t=\operatorname{Im} w - |z|^2}$$

for any $\phi \in C^\infty(0, +\infty)$.

*Proof of Theorem 3.* Taking $m = n-1$, specializing (28) further to $z = 0$ and $x = 0$, and comparing with (25), we see that

$$R(x,y) = \frac{2^{5-2n} \pi^{\frac{n-1}{2}}}{\Gamma(\frac{n-1}{2})} K(0, yi)$$

for any admissible weight function $\rho$ on $(0, +\infty)$. Note that the domain $\mathcal{S}$ is just the Cayley transform of the unit ball $\mathbf{B}^{2m}$ of $\mathbf{C}^m$, that is, although unbounded, it is biholomorphic to a bounded domain, and thus still susceptible to (1). Taking for the weight $\rho$ above the $\rho^\alpha$ from the statement of the theorem (and writing $K_\alpha$ for the corresponding $K$), and noting from the last paragraph before this proof that $\log \frac{1}{\rho(y-|z|^2)}$ is strictly-psh on $\mathcal{S}$, we thus get from (1)

$$K_\alpha(0, yi) \sim \frac{\alpha^m}{\pi^m \rho^\alpha} \det\left[\partial\overline{\partial} \log \frac{1}{\rho}\right] \qquad \text{as } \alpha \nearrow +\infty.$$

Applying (29) with $\phi = \log \frac{1}{\rho}$ therefore yields

$$K_\alpha(0, yi) \sim \frac{\alpha^{n-1}}{4\pi^{n-1}\rho^\alpha} \left(\frac{\rho'}{\rho}\right)^{n-2} \left(-\frac{\rho'}{\rho}\right)',$$

so finally

$$R_\alpha(x,y) \sim \frac{2^{3-2n}\alpha^{n-1}}{\pi^{\frac{n-1}{2}}\Gamma(\frac{n-1}{2})\rho^\alpha} \left(\frac{\rho'}{\rho}\right)^{n-2} \left(-\frac{\rho'}{\rho}\right)',$$

completing the proof of Theorem 3. □

## 5. Concluding remarks

**5.1.** With very little extra work, it can in fact be shown that both in Theorem 2 and Theorem 3, one can actually get not only the leading term but the full asymptotic expansion in decreasing powers of $\alpha$, i.e.

$$\rho(x)^\alpha R_\alpha(x) \approx \sum_{j=0}^\infty b_j(x)\alpha^{n-1-j} \qquad \text{as } \alpha \nearrow +\infty,$$

with $b_0$ given by (9) and (10). This is immediate from (2) (with $y = x$) and the last proof for Theorem 3, while for Theorem 2 the only additional item needed is that the Laplace integral (23) also admits the full asymptotic expansion

$$(30) \quad \int_a^b F(x)e^{\alpha S(x)}\,dx \approx \frac{e^{\alpha S(b)}}{\alpha S'(b)}\sum_{j=0}^\infty \Big[\Big(\frac{d}{dx}\frac{1}{S'(x)}\Big)^j F(x)\Big]_{x=b}\alpha^{-j} \qquad \text{as } \alpha \nearrow +\infty$$

if $S$ is increasing and $S'(b) > 0$ ([12], §II.1.4).

**5.2.** The case $x = 0$ omitted in Theorem 2 is easily handled directly: namely, from (17),

$$R_\alpha(0) = \frac{\Gamma(\frac{n}{2})}{2\pi^{n/2}}\Big/\int_0^1 r^{n-1}\phi(r^2)^\alpha\,dr.$$

If $\phi'(0) < 0$, then — as we know $\phi$ to be decreasing on $(0,1)$ — we again arrive at a Laplace integral that can be handled by (30), or, more precisely, by the generalization of (30) allowing $S$ with $S'(a) = 0 \neq S''(a)$, see [12], §II.1.6 (with $m = 2$ there). The result is

$$\phi(0)^\alpha R_\alpha(0) \approx \frac{\alpha^{n/2}}{\pi^{n/2}}\sum_{j=0}^\infty a_j\,\alpha^{-j/2} \qquad \text{as } \alpha \nearrow +\infty$$

with

$$a_0 = \Big(-\frac{\phi'(0)}{\phi(0)}\Big)^{n/2},$$

thus recovering, in particular, the cases $x = 0$ in (6) and (8).

Note that not only the leading power of $\alpha$ is now different than for $x \neq 0$ ($\alpha^{n/2}$ opposed to $\alpha^{n-1}$), but also the powers go down not by 1 but by $\frac{1}{2}$ in the full expansion. One can, furthermore, handle in the same way also the case when $\phi'(0) = \phi''(0) = \cdots = \phi^{(m-1)}(0) = 0 \neq \phi^{(m)}(0)$ for some $m > 1$, in which case the powers in the expansion go down by $\frac{1}{2m}$ (see again §II.1.6 in [12]). Thus the case of $x = 0$ is fundamentally different from $x \neq 0$.

Note that this jump in asymptotics (Stokes phenomenon) at $x = 0$ has no analogue on the upper half-space in Theorem 3. Apparently it seems to be connected with the fact that $\rho(x) = \phi(|x|^2)$ has a maximum at $x = 0$.

**5.3.** Using the doubling formula $\Gamma(n-1) = \pi^{-1/2}2^{n-2}\Gamma(\frac{n}{2})\Gamma(\frac{n-1}{2})$ for the Gamma function, one can rewrite the constant factor in (9) in the form

$$\frac{2\Gamma(\frac{n}{2})}{\pi^{n/2}\Gamma(n-1)} = \frac{2^{3-n}}{\pi^{\frac{n-1}{2}}\Gamma(\frac{n-1}{2})}$$

which is almost the same as the constant factor in (10), having only $2^{3-n}$ in the place of $2^{3-2n}$.

**5.4.** For $n = 2$, the right-hand sides of both (9) and (10) become simply $\frac{1}{2\pi}\Delta\log\frac{1}{\rho}$. Can it possibly be true that, for any bounded domain $\Omega \subset \mathbf{C}$ with smooth boundary and $\rho \in C^\infty(\overline{\Omega})$ such that $\rho > 0$ on $\Omega$, $\rho = 0 < \|\nabla\rho\|$ on $\partial\Omega$, and $\Delta\log\frac{1}{\rho} > 0$, one has

$$\lim_{\alpha\to+\infty} \frac{\rho^\alpha R_\alpha}{\alpha} = \frac{1}{2\pi}\Delta\log\frac{1}{\rho} \tag{31}$$

at all points where $\nabla\rho \neq 0$, while some "Stokes phenomenon" occurs at the critical points of $\rho$?

**5.5.** In some sense, the results of this paper perhaps raise more questions than they answer. On the one hand, one can conjecture, generalizing (31), that for arbitrary bounded domain $\Omega \subset \mathbf{R}^n$, $n \geq 2$, with smooth boundary, and "nice" defining function $\rho$ for $\Omega$, one has

$$\lim_{\alpha\to+\infty} \alpha^{1-n}\rho^\alpha R_\alpha = \mathcal{D}_0\rho \qquad \text{when } \nabla\rho \neq 0,$$

with some nonlinear differential operator $\mathcal{D}_0$ (depending on $\Omega$), and more generally

$$\rho^\alpha R_\alpha \approx \sum_{j=0}^\infty \alpha^{n-1-j}\mathcal{D}_j\rho \qquad \text{when } \nabla\rho \neq 0 \tag{32}$$

with some $\mathcal{D}_j$. (Thus (31) is equivalent to $\mathcal{D}_0\rho = \frac{1}{2\pi}\Delta\log\frac{1}{\rho}$ for $n = 2$.) However, it is clear neither what "nice" should mean, nor what $\mathcal{D}_0$ (or even $\mathcal{D}_j$) could look like. In fact, it is already quite surprising that $\alpha^{n-1}$ should occur as the leading power in (32), compared to the leading order $\alpha^n$ in (1) for the holomorphic case. (Note further that in (1) $n$ is the complex dimension, while in (32) $n$ is the real dimension!) Finally, one can only guess what the situation might be for other function classes in $L^2(\Omega, \rho^\alpha)$ that admit reproducing kernels, such as the pluriharmonic or the caloric functions, or quite generally the functions annihilated by a given (hypo)elliptic linear partial differential operator.

HARMONIC BERGMAN KERNELS 15

Mathematics Institute, Silesian University at Opava, Na Rybníčku 1, 74601 Opava, Czech Republic and Mathematics Institute, Žitná 25, 11567 Prague 1, Czech Republic
   *E-mail address*: englis@math.cas.cz